\documentclass[10pt, a4paper
]{amsart}

\usepackage{amsfonts, amsthm, amsmath, amscd, amssymb}
\usepackage[draft=false]{hyperref}
\usepackage[T1]{fontenc}
\usepackage{ae}
\usepackage[all]{xy}
\usepackage{xcolor}

\renewcommand{\le}{\leqslant}
\renewcommand{\ge}{\geqslant}

\newtheorem{theorem}{Theorem}

\newtheorem{proposition}[theorem]{Proposition}

\theoremstyle{definition}

\numberwithin{equation}{section}

\renewcommand\rho{\varrho}

\title{A Bombieri-Vinogradov-type theorem for moduli with small radical}


\author{Stephan~Baier}
\address{Stephan~Baier\\
	Ramakrishna Mission Vivekananda Educational Research Institute\\
	Department of Mathematics\\
	G.\ T.\ Road, PO~Belur Math, Howrah, West Bengal~711202\\
	India}
\email{stephanbaier2017@gmail.com}
\urladdr{https://www.researchgate.net/profile/Stephan\_Baier2}

\author{Sudhir Pujahari}
\address{School of Mathematical Sciences, National Institute of Science Education and Research, Bhubaneswar, An OCC of Homi Bhabha National Institute,  P. O. Jatni,  Khurda 752050, Odisha, India.}
\email{spujahari@niser.ac.in}
\urladdr{https://sites.google.com/site/sudhirkumarpujahari/home}

\begin{document}

\begin{abstract}
In this article, we extend our recent work \cite{BaPu} on a Bombieri-Vinogradov-type theorem for sparse sets of prime powers $p^N\le x^{1/4-\varepsilon}$ with $p\le (\log x)^C$ to sparse sets of moduli $s\le x^{1/3-\varepsilon}$ with radical rad$(s)\le x^{9/40}$. To derive our result, we combine our previous method with a Bombieri-Vinogradov-type theorem for general moduli $s\le x^{9/40}$ obtained by Roger Baker. 
\end{abstract}

\maketitle

\bigskip

\tableofcontents

\section{Introduction and main results} \label{results}
Let $q$ be a positive integer and $a$ be an integer coprime to $q$. By $\mathbb{P}$ denote the set of primes. For $x\ge 2$ let 
$$
\pi(x;q,a):=\sharp\{p\le x : p\in \mathbb{P} \mbox{ and } p\equiv a \bmod{q}\}
$$
and set
$$
F(x;q,a):=\pi(x;q,a)- \frac{1}{\varphi(q)} \int\limits_2^x \frac{dt}{\log t},
$$
which is the error term in the prime number theorem for the arithmetic progression $a \bmod{q}$. Further, define
$$
F(x,q):=\max\limits_{\substack{a\\ (a,q)=1}} |F(x;q,a)| \quad \mbox{and} \quad  F^{\ast}(x,q):=\max\limits_{y\le x} |F(y,q)|.
$$ 
Set $\mathcal{L}:=\log x$ throughout the sequel. The Bombieri-Vinogradov theorem implies that 
$$
F^{\ast}(x,q)\le \frac{\pi(x)}{\varphi(q)\mathcal{L}^A}
$$
for all integers $q\in (Q,2Q]$ with at most $O(Q\mathcal{L}^{-A})$ exceptions, provided that $Q\le x^{1/2}\mathcal{L}^{-2A-6}$. Under GRH, the above inequality would hold for all $q\le x^{1/2-\varepsilon}$. In slightly modified form (for $\psi(x;q,a)$ in place of $\pi(x;q,a)$), Roger Baker \cite{Bak} proved the following result, significantly restricting the set of exceptional moduli. The cost is that $Q$ is restricted to a smaller interval as well. 

\begin{theorem}[Baker] \label{Baker} Let $Q\le x^{9/40}$. Let $\mathcal{S}$ be a set of pairwise relatively prime integers in $(Q,2Q]$. Then the number of $q$ in $\mathcal{S}$ for which
$$
F^{\ast}(x,q)>\frac{\pi(x)}{\varphi(q)\mathcal{L}^A}
$$ 
is $O(\mathcal{L}^{34+A})$. 
\end{theorem}

(The above result follows easily using partial summation from Baker's original result.) 
In \cite{BaPu}, we extended the above range to $Q\le x^{1/4-\varepsilon}$ for the case when $\mathcal{S}$ consists of powers of primes $p$ such that $p\le \mathcal{L}^C$ for some fixed but arbitrary constant $C>0$. In slightly modified form (with $x$ in place of $\pi(x)$), we established the following result.

\begin{theorem}[Baier-Pujahari, 2022] \label{previousresult} Fix $C\ge 6$ and $\varepsilon>0$. Assume that $x^{\varepsilon}\le Q\le x^{1/4-\varepsilon}$. Let $\mathcal{S}\subseteq (Q,2Q]\cap \mathbb{N}$ be a set of powers of distinct primes $p\le \mathcal{L}^C$. Then the number of $q$ in $\mathcal{S}$ for which 
$$
F^{\ast}(x,q)>\frac{\pi(x)}{\varphi(q)\mathcal{L}^A}
$$
is $O_{\varepsilon}\left(\mathcal{L}^{14+2A}\right)$. 
\end{theorem}

(Here we decided to put $\pi(x)$ in place of $x$ in the numerator for an aesthetic reason: The main term in the prime number theorem for arithmetic progressions with modulus $q$ is approximately of size $\pi(x)/\varphi(q)$. Moreover, the inequalities in Theorems \ref{Baker} and \ref{previousresult} above coincide.)
In the present article, we substantially extend the family of admissible sets $\mathcal{S}$ in Theorem \ref{previousresult}. We prove the following. 

\begin{theorem} \label{mainresult}
Fix $\varepsilon>0$ and $A>0$. Assume that $x\ge 3$ and $x^{\varepsilon}\le Q\le x^{1/3}\mathcal{L}^{-15-2A}$. Let $\mathcal{S}\subseteq (Q,2Q]\cap \mathbb{N}$ be a set of relatively prime integers $s$ with radicals $\mbox{\rm rad}(s)\le x^{9/40}$ such that
$$
F^{\ast}(x,s)>\frac{\pi(x)}{\varphi(s)\mathcal{L}^A}
$$
Then 
$$
\sharp\mathcal{S}=O_{\varepsilon}\left(\mathcal{L}^{36+A}+\mathcal{L}^{14+2A}\left(1+Q^2x^{-1/2}\right)\right).
$$
\end{theorem}

Here we define the radical $\mbox{rad}(s)$ of an integer $s$ as the product of its prime divisors, i.e.
$$
\mbox{rad}(s):=\prod\limits_{p|s} p.
$$
(We reserve the symbol $p$ for primes throughout this article.) In particular, under the condition $Q=x^{1/4-\varepsilon}$ in Theorem \ref{previousresult}, we have $\sharp\mathcal{S}=O\left(\mathcal{L}^{36+A}+\mathcal{L}^{14+2A}\right)$. If $A\ge 22$, then this agrees with the bound $O\left(\mathcal{L}^{14+2A}\right)$ in Theorem \ref{previousresult} but holds under the much weaker condition rad$(s)\le x^{9/40}$ in place of $p\le \mathcal{L}^C$. 

The idea behind the proof of  the above Theorem \ref{previousresult} in \cite{BaPu} was to use Harman's sieve to compare the number of primes in the sets 
$$
\mathcal{A}:=\{n\le y \ :\ n\equiv e \bmod{p^N}\}
$$ 
and 
$$
\mathcal{B}:=\{n\le y \ :\  n\equiv d \bmod{p}\},
$$
where we assumed $2\le y\le x$, $p\le \mathcal{L}^C$ and $e\equiv d \bmod{p}$ so that $\mathcal{A} \subseteq \mathcal{B}$. In this article, we extend this idea. Here we consider residue classes 
$$
\mathcal{A}':=\{n\le y \ :\ n\equiv e \bmod{s}\}
$$ 
and 
$$
\mathcal{B}':=\{n\le y \ :\  n\equiv d \bmod{q}\},
$$
where we assume $2\le y\le x$, $q=\mbox{rad}(s)\le x^{9/40}$ and $e\equiv d\bmod{q}$ so that $\mathcal{A}\subseteq \mathcal{B}$. In \cite{BaPu}, we controlled the cardinality of primes in $\mathcal{B}$ by the Siegel-Walfisz theorem. This forced us to choose $p$ rather small, namely $p\le \mathcal{L}^C$. Here we use Baker's Theorem \ref{Baker} above to control the cardinality of primes in $\mathcal{B}'$. By the said theorem, this cardinality satisfies the predicted asymptotic for all $q$ in $\{\mbox{rad}(s): s\in \mathcal{S}\}$ with a small number of exceptions, which we discard. Our method is essentially the same as in \cite{BaPu}, where $p$ is replaced by $q$ and $p^N$ by $s$. Even though we repeat arguments from \cite{BaPu}, we shall present our proofs in full detail for self-containedness. More on the history and context of the problem can be found in \cite{BaPu} and \cite{Bak}.\\ \\
{\bf Acknowledgements:} The authors would like to thank the anonymous referee for his valuable comments, in particular, his suggestions which led to a considerably stronger main result than what we initially thought of. (In an earlier version, we considered only the range $x^{\varepsilon}\le Q\le x^{1/4-\varepsilon}$, whereas the new main result suggested by the referee handles the much larger range $x^{\varepsilon}\le Q\le x^{1/3-\varepsilon}$.)
  
\section{Notations and preliminaries}
We shall apply Theorem \ref{Baker} to generate primes in a residue class $d_q \bmod{q}$, where we assume that $q$ is square-free and satisfies $q\le x^{9/40}$.
Then for $s$ satisfying $\mbox{rad}(s)=q$, we sieve for primes in a residue class $e_q \bmod{s}$ contained in the residue class $d_q \bmod{q}$, i.e. with $e_q\equiv d_q\bmod{q}$. Here we note that it is legitimate to write $e_q$ rather than $e_s$ because there is a one to one correspondence between $s$ and $q$ due to the coprimality of distinct elements $s$ of $\mathcal{S}$ in Theorem 3.  We also recall that $\mathcal{S}\subset (Q,2Q]$ with $x^{\varepsilon}\le Q\le x^{1/3}\mathcal{L}^{-15-2A}$ and $(d_q,q)=1$ (and hence $(e_q,s)=1$). For $y_q\le x$, we set  
\begin{equation}\label{Adef}
\mathcal{A}_q:=\left\{n\le y_q: n \equiv e_q \bmod{s}\right\},
\end{equation}
and 
\begin{equation} \label{Bdef}
\mathcal{B}_q:=\left\{n\le y_q: n \equiv d_q \bmod{q}\right\}.
\end{equation}
If $\mathcal{M}$ is a finite set of integers and $z\ge 1$, we use the notation
$$
S(\mathcal{M},z):=\sharp\{n\in \mathcal{M} : p|n \mbox{ prime } \Rightarrow p\ge z\},
$$
which is common in sieve theory. We note that
\begin{equation} \label{SpiA}
\pi(y_q;s,e_q)=S(\mathcal{A}_q,x^{1/2})+O\left(\frac{x^{1/2}}{s}+1\right)
\end{equation}
and 
\begin{equation} \label{SpiB}
\pi(y_q;q,d_q)=S(\mathcal{B}_q,x^{1/2})+O\left(\frac{x^{1/2}}{q}+1\right).
\end{equation}
An application of Buchstab's identity gives 
\begin{equation} \label{Buch}
\begin{split}
S(\mathcal{C}_q,x^{1/2})= & S(\mathcal{C}_q,x^{1/3})-\sum\limits_{x^{1/3}\le p<x^{1/2}} S(\mathcal{C}_{q,p},p)\\
= & S(\mathcal{C}_q,x^{1/3})-\sum\limits_{x^{1/3}\le p<x^{1/2}} S\left(\mathcal{C}_{q,p},x^{1/3}\right)
\end{split}
\end{equation}
for $\mathcal{C}_q=\mathcal{A}_q \mbox{ or } \mathcal{B}_q$, where we write
$$
\mathcal{C}_{q,r}:=\left\{\frac{n}{r} : n\in \mathcal{C}_q, \ r|n\right\}.
$$
In particular, $\mathcal{C}_{q,1}=\mathcal{C}_q$. 
To deduce information on $S\left(\mathcal{A}_{q,r},x^{1/3}\right)$ from information on $S\left(\mathcal{B}_{q,r},x^{1/3}\right)$, we use Harman's sieve with averaging of $q$ over a subset $\mathcal{Q}$ of the set $\{\mbox{rad}(s):s\in \mathcal{S}\}$. We recall that this set consists of square-free numbers $q\le x^{9/40}$.  We apply the following extension of \cite[Proposition 7]{BaPu}.  

\begin{proposition}[Version of Harman's sieve with additional averaging]\label{Harmanasymp1}
Suppose that for any sequences $(a_m)_{m\in \mathbb{N}}$ and $(b_n)_{n\in \mathbb{N}}$ of complex numbers that satisfy $|a_m|\le 1$ and $|b_n|\le 1$ we have, for some $\lambda>0$, $\alpha>0$, $\beta\le 1/2$, $M\ge 1$ and $Y\ge 1$, that
\begin{equation} \label{typeI1}
\sum\limits_{q\in \mathcal{Q}} \left| \sum\limits_{\substack{mn\in \mathcal{A}_q\\ m\le M}} a_m -\lambda 
\sum\limits_{\substack{mn\in \mathcal{B}_q\\ m\le M}} a_m \right|^2 \le Y
\end{equation}
and
\begin{equation} \label{typeII1}
\sum\limits_{q\in \mathcal{Q}} \left| \sum\limits_{\substack{mn\in \mathcal{A}_q\\ x^{\alpha}< m\le x^{\alpha+\beta}}} a_mb_n- \lambda\sum\limits_{\substack{mn\in \mathcal{B}_q\\ x^{\alpha}< m\le x^{\alpha+\beta}}} a_mb_n \right|^2 \le Y.
\end{equation}
Let $(c_r)_{r\in \mathbb{N}}$ be a sequence of complex numbers such that 
$$
|c_r| \le 1, \mbox{ and if } c_r\not=0, \mbox{ then } p|r,\ p\in \mathbb{P} \Rightarrow p > x^{\varepsilon},
$$
for some $\varepsilon> 0$. Then, if $x^{\alpha}<M$, $2R <\min\left\{x^{1-\alpha},M\right\}$, and $M>x^{1-\alpha}$ if $2R>x^{\alpha+\beta}$, we have 
\begin{equation} \label{Sasymp1}
\sum\limits_{q\in \mathcal{Q}} \left| \sum\limits_{r\sim R} c_r\left(S\left(\mathcal{A}_{q,r},x^{\beta}\right)-\lambda S\left(\mathcal{B}_{q,r},x^{\beta}\right)\right)\right|^2 = O\left(Y\mathcal{L}^6\right).
\end{equation}
\end{proposition}

This is \cite[Theorem 3.1]{HarPri} with additional averaging over $q$ and can be established along the same lines. Here it is important that the coefficients $a_m$ and $b_n$ arising from the decompositions of $\mathcal{S}(\mathcal{C}_{q,r},x^{\beta})$ are independent of $q$. Care needs to be taken when applying \cite[Lemma 2.2]{HarPri}. Here an integral occurs which needs to be pulled out of the modulus square. This is done via an application of the Cauchy-Schwarz inequality for integrals and creates an extra factor of $\mathcal{L}$. 

Following usual custom, we call the bilinear sums in \eqref{typeI1} type I sums
and the bilinear sums in \eqref{typeII1} type II sums. As in \cite{BaPu}, we shall obtain a satisfactory type I estimate by elementary means and a satisfactory type II estimate by using a dispersion argument, followed by an application of the large sieve after detecting the implicit congruence relations using Dirichlet characters. Below is a version of the large sieve for Dirichlet characters (see \cite[Satz 5.5.1.]{Bru}, for example).

\begin{proposition}[Large Sieve] \label{ls} Let $Q$ and $N$ be positive integers and $M$ be an integer. Then, we have
$$
\sum\limits_{q\le Q} \frac{q}{\varphi(q)} \sideset{}{^\ast}\sum\limits_{\chi \bmod q} \left| \sum\limits_{M<n\le M+N} a_n\chi(n)\right|^2 \le \left(Q^2+N-1\right) \sum\limits_{M<n\le M+N} |a_n|^2,
$$
where the asterisk indicates that the sum is restricted to primitive characters. 
\end{proposition}

Another technical devise which we shall use to make ranges of variables independent is the following approximate version of Perron's formula (see \cite[Lemma 1.4.2.]{Bru}). 

\begin{proposition}[Perron's formula] \label{Perron} Let $c>0$, $N\ge 2$ and $T\ge 2$. Let $(c_n)_{n\in \mathbb{N}}$ be a sequence of complex numbers and assume that the corresponding Dirichlet series $\sum\limits_{n=1}^{\infty} c_nn^{-z}$ converges absolutely for $z=c$. Then
$$
\sum\limits_{n\le N} c_n=\frac{1}{2\pi i} \int\limits_{c-iT}^{c+iT} \left(\sum\limits_{n=1}^{\infty} c_nn^{-z}\right) N^z \frac{dz}{z}+
O\left(\frac{N^c}{T} \sum\limits_{n=1}^{\infty} |c_n|n^{-c} + C_N\left(1+\frac{N\log N}{T}\right)\right),
$$
where 
$$
C_N:=\max\limits_{3N/4\le n\le 5N/4} |c_n|. 
$$
\end{proposition}

\section{Proof of Theorem \ref{mainresult}}
We take 
$$
\mathcal{Q}:=\left\{q=\mbox{rad}(s) : s\in \mathcal{S}, \ F^{\ast}(x,q)\le \frac{\pi(x)}{\varphi(q)\mathcal{L}^{A+2}}\right\}.
$$
Using Theorem \ref{Baker}, we have 
\begin{equation} \label{exceptionsbound}
\sharp(\{\mbox{rad}(s) : s\in \mathcal{S}\}\setminus \mathcal{Q}) = O\left(\mathcal{L}^{36+A}\right). 
\end{equation}
To prove Theorem \ref{mainresult}, we apply Proposition \ref{Harmanasymp1} with 
$$
c_r=\begin{cases} 1 & \mbox{ if } r\in \{1\}\cup \{p\in \mathbb{P}: p>x^{\varepsilon}\}\\
0 & \mbox{ otherwise,}\end{cases}
$$
$$
R=1 \mbox{ or } x^{1/3}\le R<x^{1/2}, \quad \alpha=1/3=\beta, \quad  M=2x^{1/2}+1.
$$ 
Here the sets $\mathcal{A}_q$ and $\mathcal{B}_q$ are defined as in \eqref{Adef} and \eqref{Bdef}, where we choose $e_q$ and $y_q$ in such a way that
$$
|F(y_q;s,e_q)|=F^{\ast}(x,s).
$$  
We shall establish \eqref{typeI1} and \eqref{typeII1} with 
\begin{equation} \label{choice}
\lambda:=\frac{q}{s} \quad \mbox{and} \quad Y:=\left(\frac{x^2}{Q^2}+x^{3/2}+\frac{x^{5/3}\sharp\mathcal{Q}}{Q}+x\sharp\mathcal{Q}\right)\mathcal{L}^4,
\end{equation}
where we recall that $s\in (Q,2Q]$ for all $q\in \mathcal{Q}$. 
Using the inequality $|a|^2\le 2(|a-b|^2+|b|^2)$, valid for all complex $a$ and $b$, and the definitions of $F(y_q;q,d_q)$, $F(y_q;s,e_q)$ and $\mathcal{Q}$, we obtain
\begin{equation} \label{comb1}
\begin{split}
\sum\limits_{q\in \mathcal{Q}} \left|F(y_q;s,e_q)\right|^2\le  &2\sum\limits_{q\in \mathcal{Q}} \left|F(y_q;s,e_q)-\frac{q}{s}F(y_q;q,d_q)\right|^2+
2\sum\limits_{q\in \mathcal{Q}} \left|\frac{q}{s}F(y_q;q,d_q)\right|^2\\
= & 2\sum\limits_{q\in \mathcal{Q}} \left|\pi(y_q;s,e_q)-\frac{q}{s} \pi(y_q;q,d_q)\right|^2+O\left(\frac{\pi(x)^2\sharp\mathcal{Q}}{Q^2\mathcal{L}^{2B}}\right),
\end{split} 
\end{equation}
where we note that 
$$
\frac{\varphi(q)}{\varphi(s)}=\frac{q}{s}.
$$
Further, using \eqref{SpiA}, \eqref{SpiB}, \eqref{Buch} and \eqref{Sasymp1}, and using the Cauchy-Schwarz inequality after breaking up the $p$-sum into sums over dyadic intervals, we have 
\begin{equation} \label{comb2}
\sum\limits_{q\in \mathcal{Q}} \left|\pi(y_q;s,e_q)-\frac{q}{s} \pi(y_q;q,d_q)\right|^2\ll \left(\frac{x^2}{Q^2}+x^{3/2}+\frac{x^{5/3}\sharp\mathcal{Q}}{Q}+x\sharp\mathcal{Q}\right)\mathcal{L}^{12}
\end{equation}
for $Q$ as in Theorem \ref{mainresult}. Combining \eqref{comb1} and \eqref{comb2}, and using 
$$
|F(y_q;s,e_q)|> \frac{\pi(x)}{Q\mathcal{L}^A}
$$ 
by assumption in Theorem \ref{mainresult}, it follows that
$$
\frac{\pi(x)^2\sharp\mathcal{Q}}{Q^2\mathcal{L}^{2A}}\ll \sum\limits_{q\in \mathcal{Q}} \left|F(y_q;s,e_q)\right|^2\ll \left(\frac{x^2}{Q^2}+x^{3/2}+\frac{x^{5/3}\sharp\mathcal{Q}}{Q}+x\sharp\mathcal{Q}\right)\mathcal{L}^{12}+\frac{\pi(x)^2\sharp\mathcal{Q}}{Q^2\mathcal{L}^{2(A+1)}}
$$
and hence 
$$
\sharp \mathcal{Q}=O\left(\left(Q^2x^{-1/2}+1\right)\mathcal{L}^{14+2A}+\sharp\mathcal{Q}\mathcal{L}^{-1}\right)
$$  
under the condition $Q\le x^{1/3}\mathcal{L}^{-15-2A}$. 
The statement of Theorem \ref{mainresult} follows now by
discarding the elements $s$ of $\mathcal{S}$ for which $\mbox{rad}(s)\not\in \mathcal{Q}$. In view of \eqref{exceptionsbound}, the number of these exceptional $s$ is bounded by $O\left(\mathcal{L}^{36+A}\right)$.
What remains to prove are the bounds \eqref{typeI1} and \eqref{typeII1} with $Y$ defined in \eqref{choice}. This will be carried out in the remainder of this article. 
  
\subsection{Treatment of type I sums}
In the following, we suppress the index $q$ at $\mathcal{A}_q$, $\mathcal{B}_q$, $d_q$, $e_q$ and $y_q$ for simplicity. We shall treat the type I sums similarly as in \cite{BaPu}. The difference of double sums over $m$ and $n$ in \eqref{typeI1} equals
\begin{equation*} 
\Sigma_{I}=\sum\limits_{\substack{m\le M\\ (m,q)=1}} a_m\left(\sum\limits_{\substack{n\\ mn\in A}} 1 - \lambda \sum\limits_{\substack{n\\ mn\in B}} 1 \right)
\end{equation*}
which in our setting (recall the choice of $\lambda$ in \eqref{choice}) takes the form
\begin{equation*}
\Sigma_{I}=\sum\limits_{\substack{m\le M\\ (m,q)=1}} a_m\left(\sum\limits_{\substack{n\le y/m\\ n \equiv e\overline{m} \bmod{s}}} 1 - \frac{q}{s} \sum\limits_{\substack{n\le y/m\\ n \equiv d\overline{m} \bmod{q}}} 1 \right), 
\end{equation*}
where $\overline{m}$ denotes a multiplicative inverse of $m$ modulo $s$. We observe that the difference contained in the sum on the right-hand side above is $O(1)$, and hence we have 
$$
\Sigma_{I}\ll  \sum\limits_{m\le M} a_m \ll M, 
$$
which gives a bound of $O(M^2\sharp\mathcal{Q})$ for the left-hand side of \eqref{typeI1}.  Hence, \eqref{typeI1} holds with $Y$ as defined in \eqref{choice} if
\begin{equation*}
M\ll \left(\frac{x}{Q(\sharp \mathcal{Q})^{1/2}}+\frac{x^{3/4}}{(\sharp \mathcal{Q})^{1/2}}+\frac{x^{5/6}}{Q^{1/2}}+x^{1/2}\right)\mathcal{L}^2. 
\end{equation*}
Recalling our choice $M=2x^{1/2}+1$ and using $\sharp\mathcal{Q}\le \sharp\mathcal{S}\le Q\le x^{1/3}$, this is satisfied and hence \eqref{typeI1} is established.

\subsection{Treatment of type II sums.}
Our treatment of type II sums is similar as in \cite{BaPu}. Suppressing again the index $q$, the difference of double sums in \eqref{typeII1} equals  
\begin{equation} \label{introSigma} 
\Sigma_{II}:=\sum\limits_{\substack{x^{\alpha}<m\le x^{\alpha+\beta}}} a_m \left(\sum\limits_{\substack{n\\ mn\in \mathcal{A}}} b_n -  \lambda\sum\limits_{\substack{n\\ mn\in \mathcal{B}}} b_n\right).
\end{equation}
In our setting,
\begin{equation} \label{Sigmadef} 
\Sigma_{II}=\sum\limits_{\substack{x^{\alpha}<m\le x^{\alpha+\beta}\\ (m,q)=1}} a_m \left(\sum\limits_{\substack{n\le y/m\\ n\equiv e\overline{m} \bmod{s}}} b_n -  \frac{q}{s}\sum\limits_{\substack{n\le y/m\\ n\equiv d\overline{m} \bmod{q}}} b_n\right).
\end{equation}
We split $\Sigma_{II}$ into $O(\mathcal{L})$ sub-sums of the form
\begin{equation} \label{aftersplit}
\Sigma(K):=\sum\limits_{\substack{K<m\le K'\\ (m,q)=1}} a_m \left(\sum\limits_{\substack{n\le y/m\\ n\equiv e\overline{m} \bmod{s}}} b_n -  \frac{q}{s}\sum\limits_{\substack{n\le y/m\\ n\equiv d\overline{m} \bmod{q}}} b_n\right)
\end{equation}
with $x^{\alpha}\le K<K'\le 2K\le x^{\alpha+\beta}$. Throughout the following, let
$$
L:=\frac{x}{K}.
$$
To disentangle the summation variables, we apply Perron's formula, Proposition \ref{Perron}, with
$$
N:=\frac{y}{m}, \quad c:=\frac{1}{\log L}, \quad T:=L\log L
$$
and 
$$
c_n:= \begin{cases} b_n & \mbox{ if } n\le L \mbox{ and } n \equiv e \overline{m} \bmod{s} 
\mbox{ (or } n \equiv d\overline{m} \bmod{q})\\ 0 & \mbox{ otherwise} \end{cases} 
$$ 
to the inner sums over $n$ on the right-hand side of \eqref{aftersplit}. This gives
\begin{equation*}
\begin{split} 
& \sum\limits_{\substack{n\le y/m\\ n\equiv e\overline{m} \bmod{s}}} b_n -  \frac{q}{s}\sum\limits_{\substack{n\le y/m\\ n\equiv d\overline{m} \bmod{q}}} b_n\\
= & \frac{1}{2\pi i} \int\limits_{c-iT}^{c+iT} \left(\sum\limits_{\substack{n\le L\\ n\equiv e\overline{m} \bmod{s}}} b_n n^{-z} -  \frac{q}{s}\sum\limits_{\substack{n\le L\\ n\equiv d\overline{m} \bmod{q}}} b_n n^{-z} \right) \left(\frac{y}{m}\right)^{z} \frac{dz}{z} + O\left(1\right).
\end{split}
\end{equation*}
Hence, we have
\begin{equation*} 
\Sigma(K)=\frac{1}{2\pi i} \int\limits_{c-iT}^{c+iT} \Sigma(K,z)\frac{dz}{z} + O(K),
\end{equation*}
where 
\begin{equation*} 
\Sigma(K,z):=\sum\limits_{\substack{K<m\le K'\\ (m,q)=1}} a_m(z) \left(\sum\limits_{\substack{n\le L\\ n\equiv e\overline{m} \bmod{s}}} b_n(z) -  \frac{q}{s}\sum\limits_{\substack{n\le L\\ n\equiv d\overline{m} \bmod{q}}} b_n(z)\right)
\end{equation*}
with 
$$
a_m(z):=a_m\cdot \left(\frac{y}{m}\right)^z, \quad b_n(z):=b_nn^{-z}.
$$
We note that $|a_m(z)|\ll |a_m|\le 1$ if $K< m\le K'$ and $|b_n(z)|\ll |b_n|\le 1$ if $n\le L$. By the Cauchy-Schwarz inequality for integrals, we deduce that
\begin{equation} \label{aftercauchy1}
\begin{split}
|\Sigma(K)|^2 \ll & \left(\int\limits_{-T}^{T} \frac{dt}{|c+it|}\right) \cdot \int\limits_{-T}^T \frac{|\Sigma(K,c+it)|^2}{|c+it|}dt + K^2\\
\ll &
\mathcal{L} \cdot \int\limits_{-T}^T \frac{|\Sigma(K,c+it)|^2}{|c+it|}dt + K^2. 
\end{split}
\end{equation}
Set $z:=c+it$. An application of the Cauchy-Schwarz inequality for sums gives
\begin{equation} \label{CS} 
|\Sigma(K,z)|^2\le \left(\sum\limits_{\substack{K<m\le K'\\ (m,q)=1}} |a_m(z)|^2\right) \cdot \Sigma'(K,z)
\le K\Sigma'(K,z),
\end{equation}
where 
\begin{equation*}
\Sigma'(K,z):=\sum\limits_{\substack{K<m\le 2K\\ (m,q)=1}} \left|\sum\limits_{\substack{n\le L\\ n\equiv e\overline{m} \bmod{s}}} b_n(z) -  \frac{q}{s}\sum\limits_{\substack{n\le L\\ n\equiv d\overline{m} \bmod{q}}} b_n(z)\right|^2.
\end{equation*}

Now we use a dispersion argument. We multiply out the modulus square and re-arrange summations to get
\begin{equation*}
\Sigma'(K,z)=\Sigma_1(K,z)-\Sigma_2(K,z)-\Sigma_3(K,z)+\Sigma_4(K,z),
\end{equation*}
where
$$
\Sigma_1(K,z):=\sum\limits_{\substack{n_1,n_2\le L\\ n_1 \equiv n_2 \bmod{s}\\ (n_1n_2,q)=1}} b_{n_1}(z)\overline{b_{n_2}(z)} 
\sum\limits_{\substack{K<m\le 2K\\ m \equiv e\overline{n_1} \bmod{s}}} 1, 
$$
$$
\Sigma_2(K,z):=\frac{q}{s} \sum\limits_{\substack{n_1,n_2\le L\\ n_1 \equiv n_2 \bmod{q}\\ (n_1n_2,q)=1}} b_{n_1}(z)\overline{b_{n_2}(z)} 
\sum\limits_{\substack{K<m\le 2K\\ m \equiv e\overline{n_1} \bmod{s}}} 1,
$$
$$
\Sigma_3(K,z):=\frac{q}{s} \sum\limits_{\substack{n_1,n_2\le L\\ n_1 \equiv n_2 \bmod{q}\\ (n_1n_2,q)=1}} b_{n_1}(z)\overline{b_{n_2}(z)} 
\sum\limits_{\substack{K<m\le 2K\\ m \equiv e\overline{n_2} \bmod{s}}} 1,
$$
and 
$$
\Sigma_4(K,z):=\left(\frac{q}{s}\right)^2 \sum\limits_{\substack{n_1,n_2\le L\\ n_1 \equiv n_2 \bmod{q}\\ (n_1n_2,q)=1}} b_{n_1}(z)\overline{b_{n_2}(z)} 
\sum\limits_{\substack{K<m\le 2K\\ m \equiv d\overline{n_1} \bmod{q}}} 1.
$$
We see immediately that
$$
\Sigma_1(K,z)=\sum\limits_{\substack{n_1,n_2\le L\\ n_1 \equiv n_2 \bmod{s}\\ (n_1n_2,q)=1}} b_{n_1}(z)\overline{b_{n_2}(z)} \left(\frac{K}{s}+O(1)\right),
$$
$$
\Sigma_2(K,z)=\frac{q}{s} \sum\limits_{\substack{n_1,n_2\le L\\ n_1 \equiv n_2 \bmod{q}\\ (n_1n_2,q)=1}} b_{n_1}(z)\overline{b_{n_2}(z)} \left(\frac{K}{s}+O(1)\right),
$$
$$
\Sigma_3(K,z)=\frac{q}{s} \sum\limits_{\substack{n_1,n_2\le L\\ n_1 \equiv n_2 \bmod{q}\\ (n_1n_2,q)=1}} b_{n_1}(z)\overline{b_{n_2}(z)} 
\left(\frac{K}{s}+O(1)\right),
$$
and 
$$
\Sigma_4(K,z)=\left(\frac{q}{s}\right)^2 \sum\limits_{\substack{n_1,n_2\le L\\ n_1 \equiv n_2 \bmod{q}\\ (n_1n_2,q)=1}} b_{n_1}(z)\overline{b_{n_2}(z)} 
\left(\frac{K}{q}+O(1)\right)
$$
so that 
\begin{equation*}
\begin{split}
\Sigma'(K,z)=& \frac{K}{s}\sum\limits_{\substack{n_1,n_2\le L\\ n_1 \equiv n_2 \bmod{s}\\ (n_1n_2,q)=1}} b_{n_1}(z)\overline{b_{n_2}(z)} -\frac{Kq}{s^2}\sum\limits_{\substack{n_1,n_2\le L\\ n_1 \equiv n_2 \bmod{q}\\ (n_1n_2,q)=1}} b_{n_1}(z)\overline{b_{n_2}(z)}\\
&  +O\left(\frac{L^2}{s}\right),
\end{split}
\end{equation*}
provided that $s\ll L$ for all $L$ in question, i.e. for $L= x/K\ge x^{1-(\alpha+\beta)}$. This is the case since $s\in (Q,2Q]$ and 
\begin{equation} \label{somecond}
Q \le x^{1-(\alpha+\beta)}=x^{1/3}. 
\end{equation}
Now we use Dirichlet characters to detect the 
congruence relations in the above sums. We thus get 
\begin{equation*} 
\begin{split}
\Sigma'(K,z) = & \frac{K}{s} \cdot \frac{1}{\varphi(s)}\sum\limits_{\chi \bmod s} \sum\limits_{n_1,n_2\le L} b_{n_1}(z)\overline{b_{n_2}(z)} \chi(n_1)\overline{\chi}(n_2)-\\ 
&
\frac{Kq}{s^2}\cdot \frac{1}{\varphi(q)} \sum\limits_{\chi' \bmod{q}}\ \sum\limits_{n_1,n_2\le L} b_{n_1}(z)\overline{b_{n_2}(z)} \chi'(n_1)\overline{\chi'}(n_2)+O\left(\frac{L^2}{s}\right).
\end{split}
\end{equation*}
Note that 
$$
\frac{K}{s}\cdot \frac{1}{\varphi(s)}=\frac{Kq}{s^2}\cdot \frac{1}{\varphi(q)}=\frac{K}{\varphi(s^2)}
$$
and hence 
\begin{equation*}
\begin{split}
\Sigma'(K,z) = & \frac{K}{\varphi(s^2)} \sum\limits_{\chi \in \mathcal{X}(s)} \sum\limits_{n_1,n_2\le L} b_{n_1}(z)\overline{b_{n_2}(z)}\chi(n_1)\overline{\chi}(n_2) +O\left(\frac{L^2}{s}\right)\\
= & \frac{K}{\varphi(s^2)}\sum\limits_{\chi\in \mathcal{X}(s)} \left|\sum\limits_{n\le L} b_n(z) \chi(n)
\right|^2+ O\left(\frac{L^2}{s}\right),
\end{split}
\end{equation*}
where $\mathcal{X}(s)$ is the set of all Dirichlet characters modulo $s$ which are not induced by a Dirichlet character modulo $q$ (in particular, $\mathcal{X}(s)$ does not contain the principal character). We may write the above as 
\begin{equation} \label{Sigma''}
\Sigma'(K,z) = \frac{K}{\varphi(s^2)} \sum\limits_{\substack{r>q\\ q|r|s}}\ \sideset{}{^{\ast}} \sum\limits_{\chi \bmod{r}} \left|  \sum\limits_{n\le L} b_n(z) \chi(n)
\right|^2+ O\left(\frac{L^2}{s}\right),
\end{equation}
where the asterisk indicates that $\chi$ ranges over all {\it primitive} characters modulo $r$. 

Now we apply the large sieve, Proposition \ref{ls}, after re-introducing the indices $q$ and summing over $q\in \mathcal{Q}$. This gives us
\begin{equation} \label{Largesieve}
\begin{split}
\sum\limits_{q\in \mathcal{Q}}  \sum\limits_{\substack{r>q\\ q|r|s}}\frac{r}{\varphi(r)}\ \sideset{}{^{\ast}} \sum\limits_{\chi \bmod{r}}\left| \sum\limits_{n\le L} b_n(z) \chi(n)
\right|^2 \le & \sum\limits_{r\le 2Q} \frac{r}{\varphi(r)}\ \sideset{}{^{\ast}} \sum\limits_{\chi \bmod{r}}\left| \sum\limits_{n\le L} b_n(z) \chi(n)\right|^2\\
 \ll & \left(Q^2+L\right)\sum\limits_{n\le L}
|b_n|^2\\
 \ll & \left(Q^2+L\right)L,
\end{split}
\end{equation}
where we recall that $s\in (Q,2Q]$. 
Noting that 
$$
\frac{K}{\varphi(s^2)}=\frac{K}{s^2}\cdot \frac{r}{\varphi(r)}
$$ 
for $q|r|s$, and recalling that $KL=x$, we deduce that
$$
\sum\limits_{q\in \mathcal{Q}} |\Sigma_q'(K,s)|\ll Q^{-2}K^{-1}x^2+ x+Q^{-1}K^{-2}x^2\sharp\mathcal{Q}.
$$
Combining this with \eqref{aftercauchy1} and \eqref{CS}, we obtain
\begin{equation} \label{bound1}
\sum\limits_{q\in \mathcal{Q}} |\Sigma_q(K)|^2 \ll \left(Q^{-2}x^2+Kx +Q^{-1}K^{-1}x^2\sharp\mathcal{Q} + K^2\sharp\mathcal{Q}\right)\mathcal{L}^2.
\end{equation}
We get a second bound for the left-hand side by reversal of roles of the variables $m$ and $n$ in the above process, where $K$ on the right-hand side is replaced by $x/K$, i.e.
\begin{equation} \label{bound2}
\sum\limits_{q\in \mathcal{Q}} |\Sigma_q(K)|^2 \ll \left(Q^{-2}x^2+ K^{-1}x^2+ Q^{-1}Kx\sharp\mathcal{Q}+K^{-2}x^2\sharp\mathcal{Q}\right)\mathcal{L}^2. 
\end{equation}
For this to hold, the condition \eqref{somecond} needs to be replaced by
\begin{equation*}
Q \le x^{\alpha}= x^{1/3},
\end{equation*}
which is satisfied.
Using \eqref{bound1} if $K\le x^{1/2}$ and \eqref{bound2} if $K\ge x^{1/2}$, and recalling that 
$x^{\alpha}\le K\le x^{\alpha+\beta}$, we deduce that
\begin{equation*}
\sum\limits_{q\in \mathcal{Q}} |\Sigma_q(K)|^2 \ll \left(Q^{-2}x^2+ x^{3/2}+ Q^{-1}x^{2-\alpha}\sharp\mathcal{Q}
+Q^{-1}x^{1+\alpha+\beta}\sharp\mathcal{Q}+x\sharp\mathcal{Q}\right)\mathcal{L}^2.
\end{equation*}
Recalling the choice
$$
\alpha=\frac{1}{3}\quad  \mbox{and} \quad  \beta=\frac{1}{3},
$$ 
and writing $\Sigma_q=\Sigma_{II}$ with $\Sigma_{II}$ as in \eqref{Sigmadef},
we therefore obtain
\begin{equation}
\sum\limits_{q\in \mathcal{Q}} |\Sigma_q|^2\ll \left(Q^{-2}x^2+ x^{3/2}+ Q^{-1}x^{5/3}\sharp\mathcal{Q}+x\sharp\mathcal{Q}\right)\mathcal{L}^4=Y
\end{equation}
using the Cauchy-Schwarz inequality again. This establishes  \eqref{typeII1} and completes the proof of Theorem \ref{mainresult}.


\begin{thebibliography}{20}
\bibitem{BaPu} S. Baier, S. Pujahari.  
\newblock \emph{A Bombieri-Vinogradov-type theorem with prime power moduli.}
\newblock
Acta Arith. 204, No. 2, 115-140 (2022).
	
\bibitem{Bak} R.C. Baker. 
\newblock \emph{
A theorem of Bombieri-Vinogradov type with few exceptional moduli.}
\newblock
Acta Arith. 195, No. 3, 313--325 (2020).

\bibitem{Bru} J. Br\"udern.
    \newblock \emph{Einf\"uhrung in die analytische Zahlentheorie.} [\emph{Introduction to analytic number theory} (German)], Berlin: Springer-Verlag. x, 238 p. (1995).

\bibitem{HarPri}
	G.~{Harman}.
	\newblock \emph{Prime-detecting sieves}.
	\newblock Princeton, NJ: Princeton University Press (2007).


\end{thebibliography}
\end{document}